\documentclass[]{svproc}        

\usepackage{graphicx}
\usepackage[utf8]{inputenc}
\usepackage{amssymb,amsmath}
\usepackage{subcaption}
\usepackage[hidelinks]{hyperref}
\usepackage{xcolor}
\usepackage{listings}
\usepackage{environ}
\usepackage{etoolbox}

\newcommand{\R}{\mathbb{R}}
\newcommand{\Nproc}{N_P}

\newcommand{\Lap}{\nabla^2}
\newcommand{\Fprop}{\mathcal{F}}
\newcommand{\Gprop}{\mathcal{G}}
\newcommand{\NFine}{N_\Fprop}
\newcommand{\NGoarse}{N_\Gprop}

\newcounter{waycount}
\NewEnviron{way}%
{\vspace{\baselineskip}%
\par \refstepcounter{waycount}%
\parbox{\dimexpr \textwidth-1cm}%
{{\bfseries Way {\thewaycount}.} \emph{\BODY}}%
\par\addvspace{\baselineskip}}

\makeatletter
\newcommand*\NoIndentAfterEnv[1]{%
  \AfterEndEnvironment{#1}{\par\@afterindentfalse\@afterheading}}
\makeatother
\NoIndentAfterEnv{way}

\begin{document}
\mainmatter

\title{Twelve Ways To Fool The Masses When Giving Parallel-In-Time Results}
\author{Sebastian G\"otschel\inst{1} \and Michael Minion\inst{2} \and Daniel Ruprecht\inst{1} \and Robert Speck\inst{3}}

\institute{
Hamburg University of Technology\\
Institute of Mathematics\\
Chair Computational Mathematics\\
Am Schwarzenberg-Campus 3\\
D-21073 Hamburg, Germany\\
\email{\{ruprecht, sebastian.goetschel\}@tuhh.de}
\and
Lawrence Berkeley National Laboratory\\
1 Cyclotron Road\\
Berkeley, CA 94720, USA\\
\email{mlminion@lbl.gov}
\and
Forschungszentrum Jülich GmbH \\
Jülich Supercomputing Centre \\
52425 Jülich, Germany \\
\email{r.speck@fz-juelich.de}
}

\date{September 2020}

\maketitle

\setcounter{footnote}{0} 

\begin{abstract}
Getting good speedup---let alone high parallel efficiency---for parallel-in-time (PinT) integration examples can be frustratingly difficult. The high complexity and large number of parameters in PinT methods can easily (and unintentionally) lead to numerical experiments that overestimate the algorithm's performance. In the tradition of Bailey's article ``Twelve ways to fool the masses when giving performance results on parallel computers", we discuss and demonstrate pitfalls to avoid when evaluating performance of PinT methods. Despite being written in a light-hearted tone, this paper is intended  to raise awareness that there are many ways to unintentionally fool yourself and others and that by avoiding these fallacies more meaningful PinT performance results can be obtained.
\end{abstract}

\section{Introduction}
The trend towards extreme parallelism in high-performance computing requires 
novel numerical algorithms to translate the raw computing power of hardware into application performance~\cite{DongarraEtAl2014}.
Methods for the approximation of time-depen\-dent partial differential equations, which are used in models in a very wide range of disciplines from engineering to physics, biology or even sociology, pose a particular challenge in this respect. Parallelization of algorithms discretizing the spatial dimension via a form of domain decomposition is quite natural and has been an active research topic for decades.
Exploiting parallelization in the time direction
is less intuitive as time has a clear direction of information transport. 
Traditional algorithms for temporal integration  employ a step-by-step procedure that is difficult to parallelize. In many applications, this sequential treatment of 
temporal integration has become a bottleneck in massively parallel simulations.

Parallel-in-time (PinT) methods, i.e., methods that offer at least some degree of concurrency, are advertised as a possible solution to this temporal bottleneck.
The concept was pioneered by Nievergelt in 1964~\cite{Nievergelt1964}, but has only really gained traction in the last two decades~\cite{Gander2015_Review}.
By now, the effectiveness of PinT has been well established 
for  examples ranging from the linear heat equation in one-dimension to more complex highly diffusive problems in more than one dimension.
More importantly, there is now ample evidence that different PinT methods can deliver measureable reduction in solution times on real-life HPC systems for a wide variety of problems. 
Ong and Schroder~\cite{ong2019applications} and Gander~\cite{Gander2015_Review} provide overviews of the literature, and a good resource for further reading is also given by the community website \url{https://parallel-in-time.org/}.

PinT methods differ from space-parallel algorithms or parallel methods for operations like
the FFT in that they do not simply parallelize a serial algorithm to reduce its run time.\footnote{One exception are so-called ``parallel-across-the-method'' PinT methods in the terminology by Gear~\cite{Gear1988} that can deliver smaller-scale parallelism.}
Instead, serial time-stepping is usually replaced with a computationally more costly 
and typically iterative  procedure that is amenable to parallelization.
Such a procedure will run much slower in serial, but can overtake serial time-stepping in speed if sufficiently many processors are employed.
This makes a fair assessment of performance much harder since there is no clear baseline to compare against.
Together with the large number of parameters and inherent complexities in PinT methods and 
PDEs themselves, there are thus many sometimes subtle ways to fool oneself (and the masses) when assessing performance.
We will demonstrate various ways to produce results that seem to demonstrate speedup but are essentially meaningless.
The paper is written in a similar spirit as other ``ways to fool the masses" papers 
first introduced in~\cite{Bailey1991-pr} who inspired a series of similarly helpful papers in related areas~\cite{Gustafson1991-wi,Hoefler2018-wb,Minhas2019-rm,299418,Tautges2004-ry,Chawner2011-bj,Pakin2011-ec}.  One departure from the canon here is that we provide actual examples to demonstrate the {\it Ways} as we present them. Despite the light-hearted, sometimes even sarcastic tone of the discussion, the numerical examples are similar to experiments one could do for evaluating the performance of PinT methods.

Some of the ways we present are specific to PinT while others, although formulated in ``PinT language'' correspond to broader concepts from parallel computing.
This illustrates another important fact about PinT: while the algorithms often dig deeply into the applied mathematics toolkit, their relevance is exclusively due to the architectural specifics of modern high-performance computing systems.
This inherent cross-disciplinarity is another complicating factor when trying to do fair performance assessments. Lastly we note that this paper was first presented in a shorter form as a conference talk at the 9th Parallel in Time Workshop held (virtually) in June, 2020.  Hence some of the {\it Ways} are  more  relevant to  live presentations,  although all should be considered in both written and live scenarios. 
In the next section we present the 12 {\it Ways} with a series of numerical examples before concluding with some more serious comments in Section 3.

\section{Fool the masses (and yourself)}

\subsection{Choose Your Problem Wisely!}\label{sec:problem} 

If you really want to impress the masses with your PinT results, you will want to show as big a parallel speedup as possible, hence you will want to use a lot of processors in the time direction.  If you are using, for example, Parareal~\cite{LionsEtAl2001}, a theoretical model for speedup is given by the expression
\begin{equation} \label{eq:pararealS}
    S_\text{theory} = \frac{\Nproc}{\Nproc \alpha + K(1+\alpha)},    
\end{equation}
where $\Nproc$ is the number of processors,  $\alpha$ is the ratio of the cost of the coarse propagator $\Gprop$  compared to the fine propagator $\Fprop$, and $K$ is the total number of iterations needed for convergence.  Hence to get a large speedup that will impress the masses, we need to choose $\Nproc$ to be large, $\alpha$ to be small, and hope $K$ is small as well.  A common choice for parareal is to 
have $\Gprop$ be one step of some method and $\Fprop$ be $N_\Fprop$  steps of the same method so that $\alpha=1/\NFine$ is small.  
But note that this already means that the total number of time steps corresponding to the serial method is now $\Nproc \NFine$. 
Hence we want to choose an equation and problem parameters for which very many time steps can be employed, while still showing good speedup without raising any suspicions that the problem is too ``easy".  The first example suggests some {\it Ways} to pull off this perilous balancing act.

In this example, we use the following nonlinear advection-dif\-fu\-sion-reaction equation 
\begin{align*}
    u_t = v u_x + \gamma uu_x+ \nu u_{xx} + \beta u (a-u)(b-u),
\end{align*}
where the constants $v, \gamma, \nu, \beta, a$, and $b$ determine the strength of each term.  In order to squeeze in the massive number of time steps we need for good speedup, we choose a long time interval over which to integrate, $t \in [0,T_F]$,
with $T_F=30$. The initial condition is given on the interval $[0,2\pi]$ by
\begin{align*}
    u(x,0) = 1 - d(1-e^{-(x-\pi)^4/\sigma}).
\end{align*}
If you are presenting this example in front of the an audience, try to get all the equations with parameters on one slide and then move on before defining them.
\begin{way}
Choose a seemingly complicated problem with lots of parameters which you define later (or not at all).
\end{way}

For the first numerical test, we choose $\Nproc=200$ processors, and use a fourth-order IMEX Runge-Kutta method and a pseudo-spectral discretization in space using 128 grid points, where the linear advection and diffusion terms are treated implicitly. We use one time step for $\Gprop$ and $\NFine=64$ steps for $\Fprop$. 
Since the method is spectrally accurate in space, it gives us cover to use a lot of time steps (more on that later).  We set the stopping criterion for Parareal to be when the increment in the iteration is below $10^{-9}$,
and $v=-0.5$, $\gamma=0.25$, $\nu=0.01$, $\beta=-5$, $a=1$, $b=0$, and $d=0.55$ (see also Appendix 1). For these values, Parareal converges on the entire time interval in 3 iterations.  The theoretical speedup given by Eq.~\ref{eq:pararealS} is 32.4.  Not bad!

\begin{figure}
     \centering
         \includegraphics[width=0.7\textwidth]{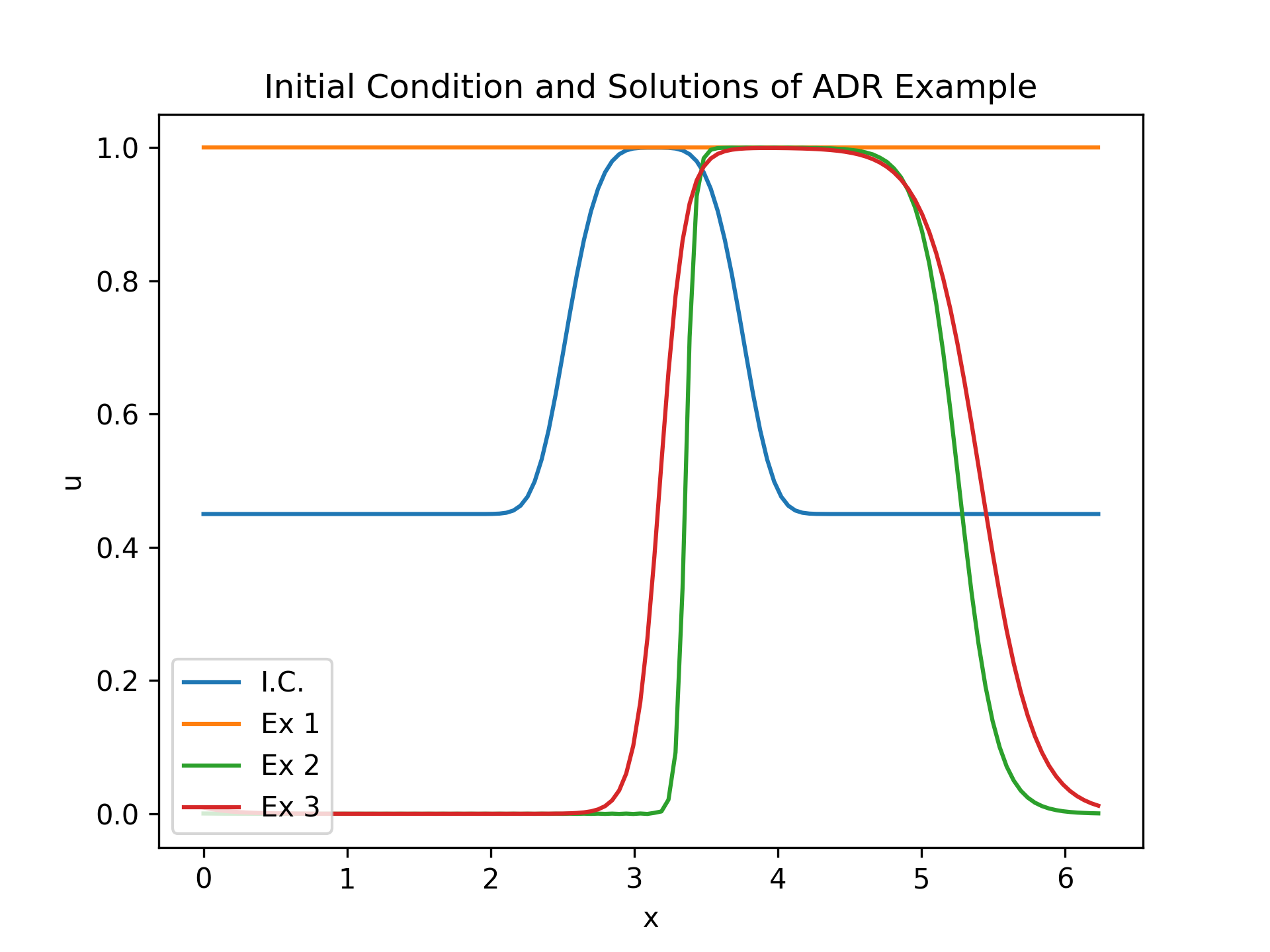}
         \caption{Initial solution and solution at $t=30$ for the advection-diffusion-reaction problem demonstrating the significant effect that parameter selection can have on the dynamics and subsequent PinT speedup discussed in {\it Ways} 1-4. }
    \label{fig:steady1}
\end{figure}

If we explore no further, we might have fooled the masses.  How did we manage?  Consider a plot of the initial condition and solution at the final time for this problem shown in Fig.~\ref{fig:steady1}, with the blue and orange lines respectively.  The lesson here is
\begin{way} \label{way:steadystate}
Quietly  use an initial condition and/or problem parameters for which the solution tends to a steady state. But do not show the actual solution.
\end{way}

If we repeat this experiment changing only one parameter, $b=0.5$, the number of Parareal iterations needed for convergence jumps to $K=10$ for a less impressive theoretical speedup of $15.05$.  In this case the solution quickly evolves not to constant state, but a steady bump moving at constant speed
(the green line in Fig. \ref{fig:steady1}).  This raises another important point to fool the masses:
\begin{way}
Do not show the sensitivity of your results to problem parameter changes.  Find the best one and let the audience think the behavior is generic.
\end{way}
Sometimes you might be faced with a situation like the second case above and not know how to get a better speedup. One suggestion is to add more diffusion.  Using the same parameters except increasing the diffusive coefficient to $\nu=0.04$ reduces the number of iterations for convergence to $K=5$ with a theoretical speedup of $24.38$.  The solution of this third example is shown by the red line
Fig. \ref{fig:steady1}.
If you can't add diffusion directly, using a diffusive discretization for advection like  first-order upwind finite differences can sometimes do the trick while avoiding the need to explicitly admit to the audience that you needed to increase the amount of ``diffusion".
\begin{way}
If you are not completely thrilled about the speedup because the number of iterations $K$ is too high, try adding more diffusion.  You might have to experiment a little to find just the right amount.
\end{way}\noindent

\subsection{Over-resolve the solution! Then over-resolve some more.}\label{sec:overresolve} 
After carefully choosing your test problem, there are ample additional opportunities to boost the parallel performance of your numerical results. The next set of {\it Ways}  consider the effect of spatial and temporal resolution.
We consider the 1D nonlinear Schrödinger equation
\begin{align}\label{eq:schroedinger}
    u_t = \mathrm{i}\Delta u + 2\mathrm{i}\left|u\right|^2u
\end{align}
with periodic boundary conditions in $[0,2\pi]$ and the exact solution as given by Aktosun et al.~\cite{aktosun2007exact}, which we also use for the initial condition at time $t_0 = 0$. 

\begin{figure}
     \centering
     \begin{subfigure}[]{0.475\textwidth}
         \centering
         \includegraphics[width=\textwidth]{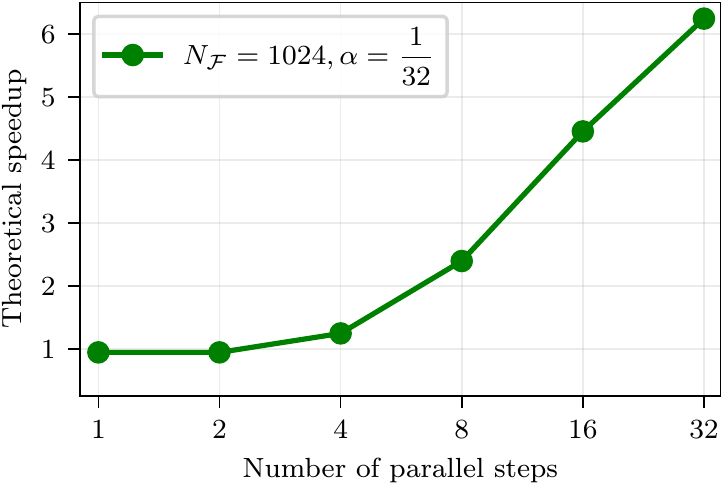}
         \caption{Deceivingly good speedup}
         \label{fig:overres_time_no}
     \end{subfigure}
     \hfill
     \begin{subfigure}[]{0.475\textwidth}
         \centering
         \includegraphics[width=\textwidth]{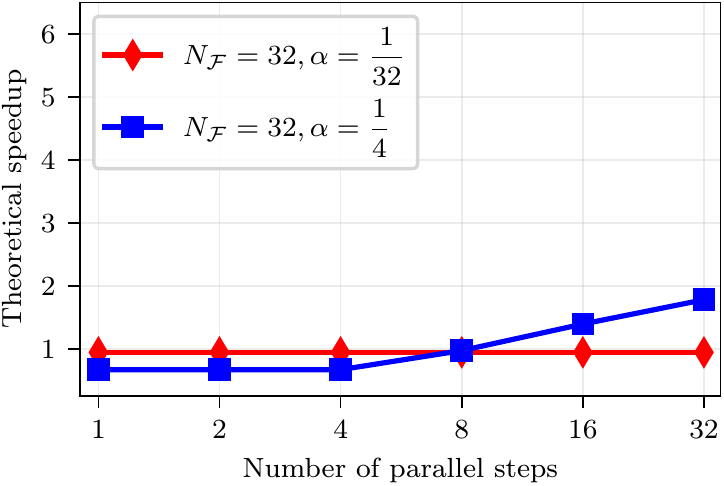}
         \caption{Not so good speedup}
         \label{fig:overres_time_yes}
     \end{subfigure}
    \caption{Estimated speedup for Parareal runs of the nonlinear Schrödinger example~\eqref{eq:schroedinger} demonstrating the effect
    of over-resolution in time (\textit{Way} 5).}
    \label{fig:overres_time}
\end{figure}

This is a notoriously hard problem for out-of-the-box PinT methods, but we are optimistic and give Parareal a try.
We use a second-order IMEX Runge-Kutta method by Ascher et al.~\cite{Ascher1997-er}
with $\NFine=1\,024$  steps and $\NGoarse=32$ coarse steps for each of the $32$ processors.  
In space, we again use a pseudo-spectral method with the linear part treated implicitly and $N_x = 32$ degrees-of-freedom.
The estimated speedup can be found in Fig.~\ref{fig:overres_time_no}. 
Using $K=5$ iterations, we obtain a solution about $6.24$ times faster when running with $32$ instead of $1$ processor. 
All runs achieve the same accuracy of $5.8\times 10^{-5}$ and it looks like speedup  in time can be easily achieved after all.  

Yet, although the accuracy compared to the exact solution is the same for all runs, 
the temporal resolution is way too high for this problem, masking the effect of 
coarsening in time.
The spatial error is dominating and instead of $32\times 1\,024 = 32\,768$ overall time steps, only $32\times 32 = 1\,024$ are actually needed to balance spatial and temporal error.
 Therefore, the coarse level already solves the problem quite well: speedup only comes from over-resolving in time.
 
If we instead choose $\NFine=32$ time steps on the fine level instead with the same coarsening factor of $32$ ($\NGoarse=1$), we get no speedup at all -- see the red curve/diamond markers in~Fig.~\ref{fig:overres_time_yes}.
Using a less drastic coarsening factor of $4$ leads to a maximum speedup of $1.78$ with $32$ processors (blue curve/square markers), which is underwhelming and frustrating and not what we would prefer to present in public.
Lesson learned:
\begin{way}
Make $\Delta t$ so small that the coarse integrator is already accurate. Never check if a larger $\Delta t$ might give you the same solution.
\end{way}
The astute readers may have noticed we also used this trick to a lesser extent in the advection-diffusion-reaction example above.

\begin{figure}
     \centering
     \begin{subfigure}[]{0.475\textwidth}
         \centering
         \includegraphics[width=\textwidth]{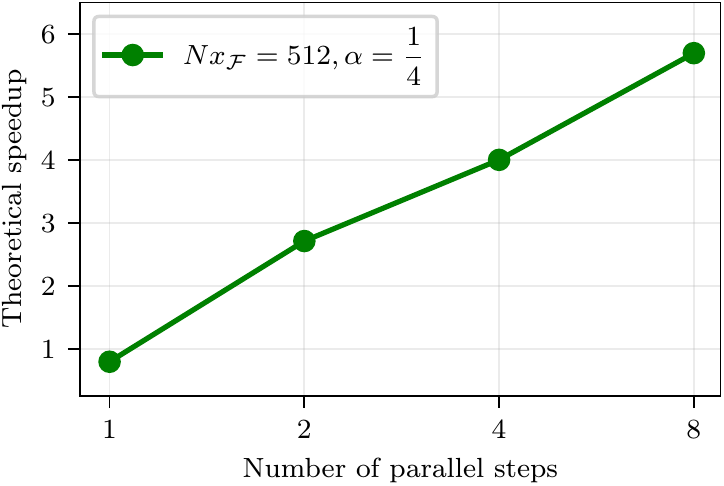}
         \caption{Deceivingly good speedup}
         \label{fig:overres_space_no}
     \end{subfigure}
     \hfill
     \begin{subfigure}[]{0.475\textwidth}
         \centering
         \includegraphics[width=\textwidth]{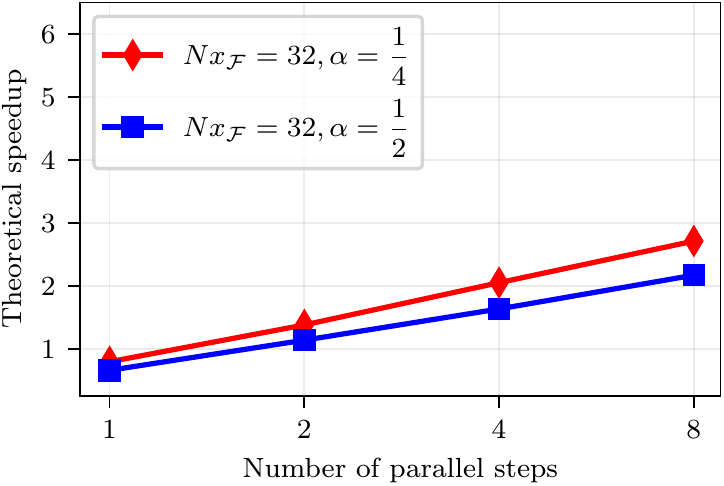}
         \caption{Not so good speedup}
         \label{fig:overres_space_yes}
     \end{subfigure}
    \caption{Estimated speedup for PFASST runs of the nonlinear Schrödinger example~\eqref{eq:schroedinger} demonstrating the effect
    of over-resolution in space ({\it Way} 6).}
    \label{fig:overres_space}
\end{figure}

A similar effect can be achieved when considering coarsening in space.
Since we have learned that more parameters are always good to fool the masses 
we now use PFASST~\cite{EmmettMinion2012} instead of Parareal.
We choose $5$ Gauss-Lobatto quadrature nodes per time step, leading to an $8$th order IMEX method, which requires only $8$ time steps to achieve an accuracy of about $5.8\times 10^{-5}$. 
We do not coarsen in time, but -- impressing everybody with how resilient PFASST is to spatial coarsening -- go from $512$ degrees-of-freedom on the fine to $32$ on the coarse level.
We are rewarded with the impressive speedup shown in Fig.~\ref{fig:overres_space_no}: using $8$ processors, we are $5.7$ times faster than the sequential SDC run.
To really drive home the point how amazing this is, we point out to the reader that this corresponds to 71\% parallel efficiency.
Even the space-parallel linear solver crowd would have to grudgingly accept such an efficiency as respectable.

However, since the spatial method did not change from the example before, we already know that $32$ degrees-of-freedom would have been sufficient to achieve a PDE error of about $10^{-5}$. So using $512$ degrees-of-freedom on the fine level heavily over-resolves the problem in space.
Using only the required $32$ degrees-of-freedom on the fine level with a similar coarsening factor of $4$ only gives a speedup of $2.7$, see Fig.~\ref{fig:overres_time_yes} (red curve/diamond markers).
While we could probably sneak this into an article, the parallel efficiency of 34\% will hardly impress anybody outside of the PinT  community.

It is worth noting that better resolution in space on the coarse level does not help (blue curve/square markers).
This is because the coarse level does not contribute anything to the convergence of the method anymore.
Turning it off completely would even increase the theoretical speedup to about $3.5$. Hence, for maximum effect:
\begin{way}
When coarsening in space, make $\Delta x$ on the fine level so small that even after coarsening, the coarse integrator is accurate.
Avoid the temptation to explore a more reasonable resolution.
\end{way}

\subsection{Be smart when setting your iteration tolerance!}\label{iterationtol}
If the audience catches on about your $\Delta t$/$\Delta x$ over-resolution issues, there is a more subtle way to over-resolve and fool the masses. Since methods like Parareal and PFASST are iterative methods, one must decide when to stop iterating - use this to your advantage!
The standard approach is to check the increment between two iterations or some sort of residual (if you can, use the latter: it sounds fancier and people will ask fewer questions).
In the runs shown above, Parareal is stopped when the difference between two iterates is below $10^{-10}$ and PFASST is stopped when the residual of the local collocation problems is below $10^{-10}$.

\begin{figure}
     \centering
     \begin{subfigure}[]{0.475\textwidth}
         \centering
         \includegraphics[width=\textwidth]{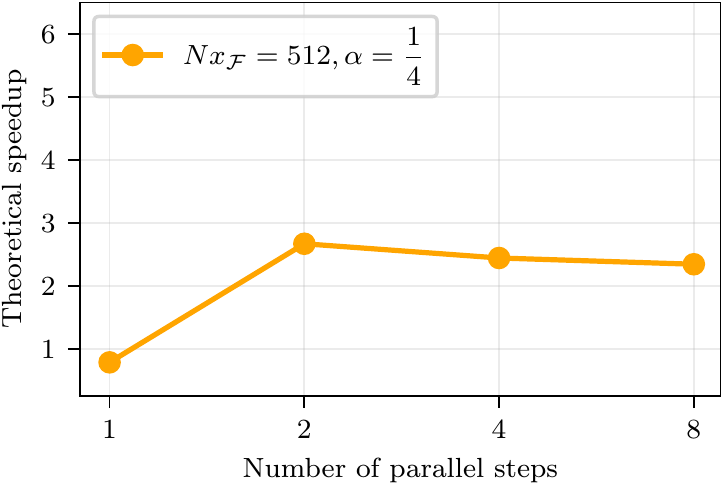}
         \caption{Now also not so good speedup}
         \label{fig:overres_iter_no}
     \end{subfigure}
     \hfill
     \begin{subfigure}[]{0.475\textwidth}
         \centering
         \includegraphics[width=\textwidth]{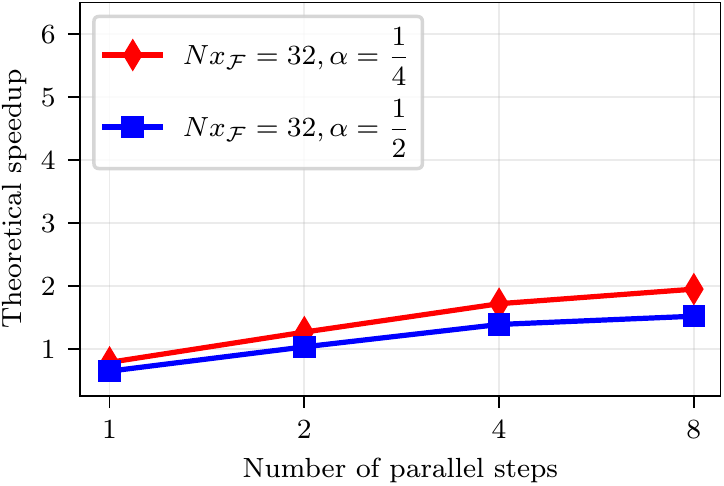}
         \caption{Still not so good speedup}
         \label{fig:overres_iter_yes}
     \end{subfigure}
    \caption{Estimated speedup for PFASST runs of the nonlinear Schrödinger example~\eqref{eq:schroedinger} with different resolutions in space demonstrating how using a sensible iteration tolerance of $10^{-5}$ can reduce speedup ({\it Way} 7).}
    \label{fig:overres_iter}
\end{figure}

These are good choices, as they give you good speedup: for the PFASST example, a threshold of $10^{-5}$ would have been sufficient to reach the accuracy of the serial method.
While this leads to fewer PFASST iterations (good!), unfortunately it also makes the serial SDC baseline much faster (bad!).
Therefore, with the higher tolerance, speedup looks much less attractive, even in the over-resolved case, see Fig.~\ref{fig:overres_iter_no}.

Similarly, when using more reasonable tolerances, the speedup of the well-resolved 
examples decreases as shown in Fig.~\ref{fig:overres_iter_yes}.  
This leads to our next {\it Way}, which has a long and proud tradition, and for which we can therefore quote Pakin~\cite{Pakin2011-ec} directly,
\begin{way}\label{way:toomanyiter}
``Hence, to demonstrate good [...] performance, always run far more iterations than are typical, necessary, practical, or even meaningful for real-world usage, numerics be damned!"
\end{way}

Yet another smart way to over-resolve  is to choose a ridiculously small tolerance for an inner iterative solver. 
Using something cool like GMRES to solve the linear systems for an implicit integrator far too accurately is an excellent 
avenue for making your serial baseline method much slower than it needs to be.
This is further desirable because it reduces the impact of tedious overheads like communication or waiting times.
We all know that an exceptional way to good parallel performance is using a really slow serial baseline to compare against.

\begin{way}
    Not only use too many outer iterations, but try to maximize the amount of work done by iterative spatial solvers (if you have one, and you always should).
\end{way}

Note that for all the examples presented so far, we did not report any actual speedups measured on parallel computers. 
Parallel programming is tedious and difficult, as everybody understands,  and 
what do we have performance models for, anyway? 
It is easier to just plug your parameters into a theoretical model.
Realizing this on an actual system can be rightfully considered Somebody Else's Problem (SEP) or a task for your dear future self.
But for completeness, the next example will address this directly.


\subsection{Don't report runtimes!} 

Because solving PDEs only once can bore an audience, we will now talk about optimal control of the heat equation, the ``hello world'' example in optimization with time-dependent PDEs.
This problems has the additional advantage that  even more parameters are available to tune.
Our specific problem is as follows. Given some desired state $u_d$ on a space-time domain $\Omega \times (0,T)$, $\Omega \subset \R^d$, we aim to find a control $c$ to minimize the objective functional
\begin{equation*}
J(u,c) = \frac{1}{2} \int_0^T \|u-u_d\|^2_{L^2(\Omega)} \ \mathrm{d}t + \frac{\lambda}{2} \int_0^T \|c\|^2_{L^2(\Omega)} \ \mathrm{d}t 
\end{equation*}
subject to \[u_t - \Lap u = c + f(u)\] 
with periodic boundary conditions 
(allowing us to use FFT to evaluate the Laplacian and perform the implicit linear solves ({\it Way} 8 be damned). 
For the linear heat equation considered in the following, the source term is $f(u) \equiv 0$ (see {\it Way} 1). 
Optimization is performed using steepest descent; for computation of the required reduced gradient we need to solve a forward-backward system of equations for state $u$ and adjoint $p$,
\begin{alignat*}{3}
    u_t - \Lap u &= c +f(u) &\qquad -p_t - \Lap p - f'(u)p &= u-u_d \\
	u(\cdot,0) &= 0  &p(\cdot,T) &= 0.
\end{alignat*}

To parallelize in time, we use, for illustration, the most simple approach: given a control $c$, the state equation is solved parallel-in-time for $u$, followed by solving the adjoint equation parallel-in-time for $p$ with PFASST using $\Nproc=20$ processors. 
For discretization, we use 20 time steps and three levels with 2/3/5 Lobatto IIIA nodes in time as well as 16/32/64 degrees of freedom in space. 
As a sequential reference we use MLSDC on the same discretization. 
We let PFASST/MLSDC iterate until the residual is below $10^{-4}$ instead of iterating to high precision, so we can 
openly boast how we avoid \textit{Way~}\ref{way:toomanyiter}.

In the numerical experiments we perform one iteration of an iterative, gradi\-ent-based optimization method to evaluate the method (i.e., solve state, adjoint, evaluate objective, compute gradient). 
As initial guess for the control we do not use the usual choice $c_0 \equiv 0$ as this would lead to $u \equiv 0$ but a nonzero initial guess---again, we make sure everybody knows that we avoid \textit{Way}~\ref{way:steadystate} in doing so.
To estimate speedup we count PFASST/MLSDC iterations and compute
\[
    S = \frac{\text{total MLSDC iterations state} + \text{adjoint}}{\text{total iterations state on } CPU_{20} + \text{total iterations adjoint on } CPU_{1}}.
\]
We get $S = \frac{40+60}{7+7}= 7.1,$ for a nice parallel efficiency of $35.5\%$.

Before we publish this, we might consider actual timings from a real computer. 
Unfortunately, using wall clock times instead of iterations gives
\[
S = \frac{\text{serial wall clock time}}{\text{parallel wall clock time}} = \frac{44.3s}{18.3s} = 2.4,
\]
and thus only roughly a third of the theoretical speedup. To avoid this embarrassment:

\begin{way}\label{way:donttime}
Only show theoretical/projected speedup numbers (but mention this only in passing or not at all). If you include the cost of communication in the theoretical model, assume it is small enough not to affect your speedup.
\end{way}

Why is the theoretical model poor here? One cause is the overhead for the optimization---after all, there is the evaluation of the objective functional, and the construction of gradient. 
Ignoring parts of your code to report better results is another proud tradition of parallel computing, see Bailey's paper~\cite{Bailey1991-pr}.
However, most of the tasks listed do trivially parallelize in time. 
The real problem is that communication on an actual HPC system is aggravatingly not really instantaneous.



\begin{figure}
    \centering
    \includegraphics[width=\textwidth]{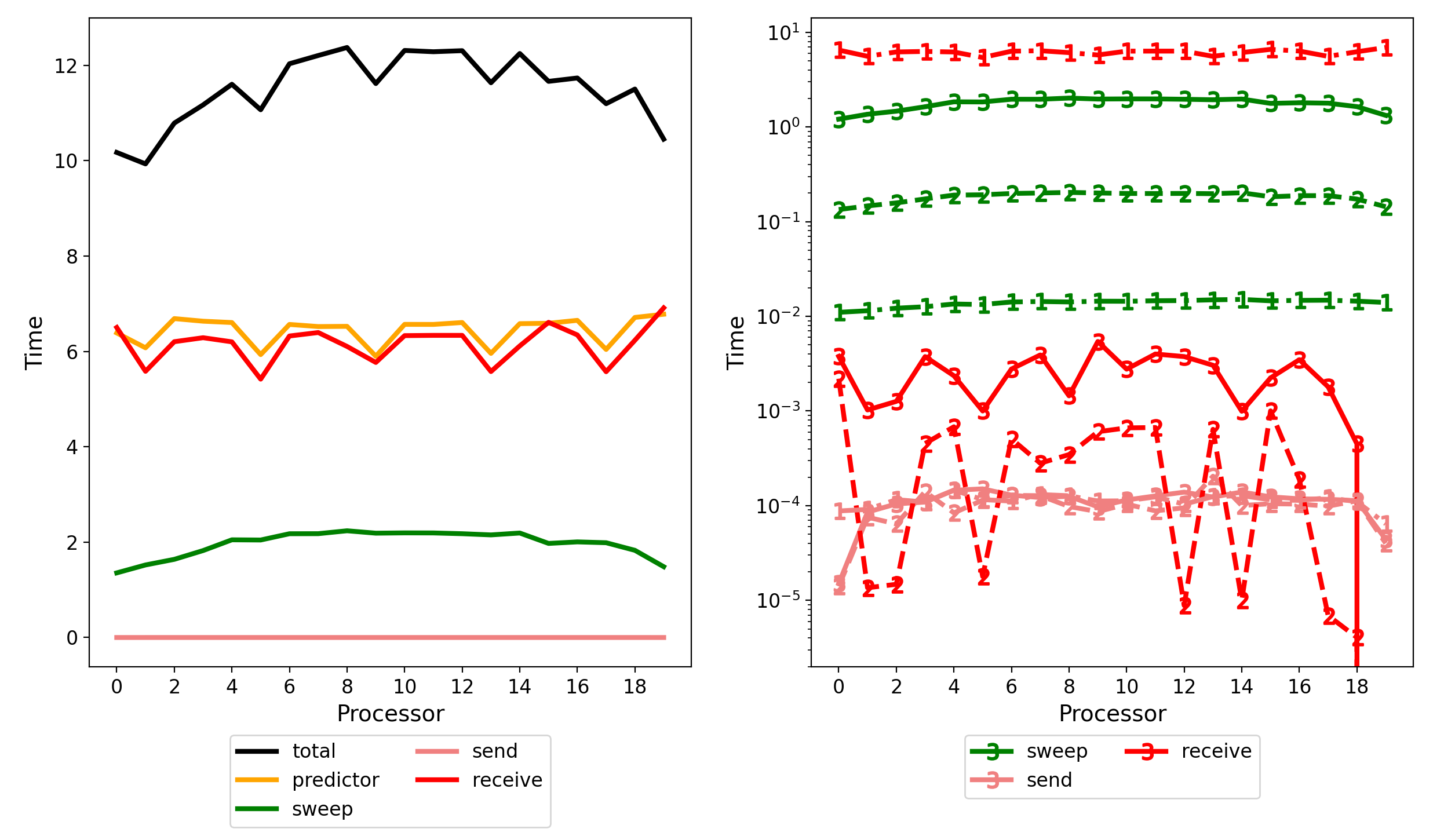}    
    \caption{Wall clock times of the different algorithmic steps for the linear heat equation example on $\Omega = [0,1]^3$ and $T=1$. Left: total times. Right: times per level (1 is coarsest level, 3 finest).  Note the "receive" times are not negligible as discussed in \textit{Way}~\ref{way:donttime}.}
    \label{fig:timings}
\end{figure}

Looking at detailed timings for PFASST, 
Fig.~\ref{fig:timings} shows that the issue truly is in communication costs, which clearly cannot be neglected. 
In fact, more time is spent on blocking coarse grid communication than on fine sweeps. 
Note also that, due to the coupled forward-backward solves, each processor requires similar computation and communication times. 
The following performance model
\[
    S = \frac{\Nproc}{\frac{\Nproc \alpha}{K_S}+ \frac{K_P}{K_S}(1+\alpha+\beta)}
\]
accounts for overheads in the $\beta$ term.
Matching the measured speedups requires setting $\beta = 3$ or three times the cost of one sweep on the fine level!
This is neither small nor negligible by any measure.


\subsection{Choose the measure of speedup to your advantage!}
Technically, parallel speedup should be defined as the ratio of the run time of your parallel code to the run time required by the \textbf{best} available serial method. But who has the time or energy for such a careful comparison?
Instead, it is convenient to choose a baseline to get as much speedup as possible.

In the example above, MLSDC was used as a baseline since it is essentially the sequential version of PFASST and allows for a straightforward comparison and the use 
of a theoretical speedup.
However, MLSDC might not be the fastest serial method to solve state and adjoint equations to some prescribed tolerance. 
For illustration, we consider solving an optimal control problem for a nonlinear heat equation with $f(u) = -\frac{1}{3}u^3 + u$ on $\Omega\times (0,T) = [0,20]\times (0,5)$. 
Wall clock times were measured for IMEX-Euler, a 4th-order additive Runge-Kutta scheme (ARK-4), and 3-level IMEX-MLSDC with 3/5/9 Lobatto IIIA collocation points, with each method reaching a similar final accuracy in the computed control (thus, using different number of time steps, but the same spatial discretization). 
For the IMEX methods, the nonlinearity as well as the control were treated explicitly. 
IMEX-Euler was fastest with $102.5s$, clearly beating MLSDC ($169.8s$) despite using significantly more time steps. 
The ARK-4 method here required $183.0s$, as the non-symmetric stage values slow down the forward-backward solves due to the required dense output.
With PFASST on $32$ CPUs requiring $32s$, the speedup reduces from $5.3$ for MLSDC as a reference to $3.2$ when compared to IMEX-Euler. 
By choosing the sequential baseline method wisely, we can increase the reported speedup in this example by more than $65\%$!

A very similar slight of hand is hidden in Sect.~\ref{sec:overresolve}, where only theoretical speedups are reported.
In the PFASST examples, the SDC iteration counts are used as the baseline results, although in most cases MLSDC required up to 50\% fewer iterations to converge. 
Using MLSDC as a baseline here would reduce the theoretical speedups significantly in all cases. 
Whether this still holds when actual run times are considered, is, as we have just seen, part of a different story.
\begin{way}
If you report speedup based on actual timings, compare your code to the method run on one processor and never against a different and possibly more efficient serial method.
\end{way}

\subsection{Use low order serial methods!}
A low-order  temporal method 
is a choice convenient for PinT methods because they are easier to implement 
and allow one to take many time steps without falling prey to  \textit{Way} 5, especially when you want to show how the speedup increases as
you take ever more time steps for a problem on a fixed time interval.  After all, it is the parallel scaling that is exhilirating, not necessarily how quickly one can compute a solution to a given accuracy.

\begin{figure}
     \centering
     \hspace{-1cm}
     \begin{subfigure}[]{0.475\textwidth} 
         \centering
         \includegraphics[width=1.2\textwidth]{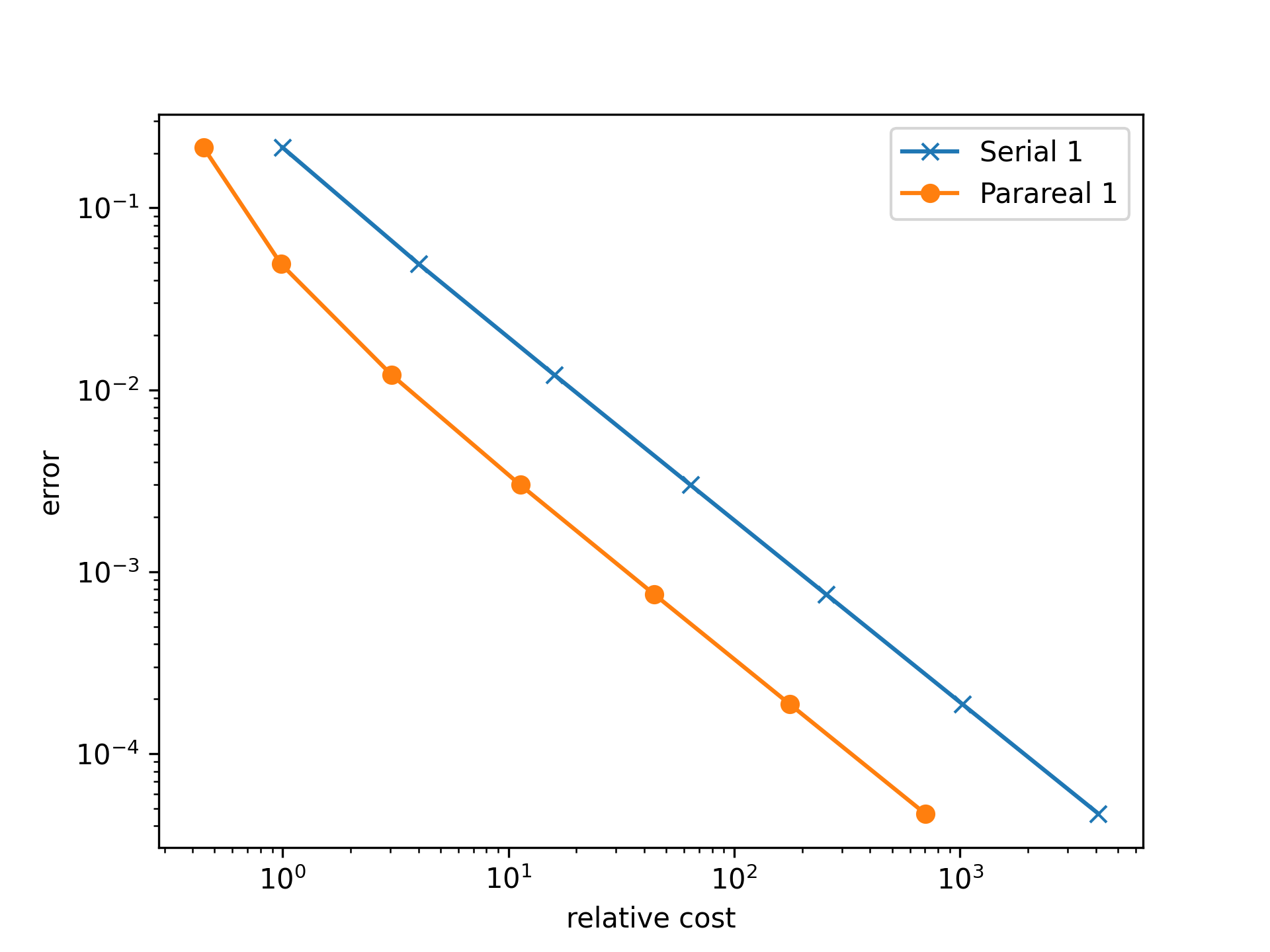}
     \end{subfigure}
     \hfill
     \begin{subfigure}[]{0.475\textwidth}
         \centering
         \includegraphics[width=1.2\textwidth]{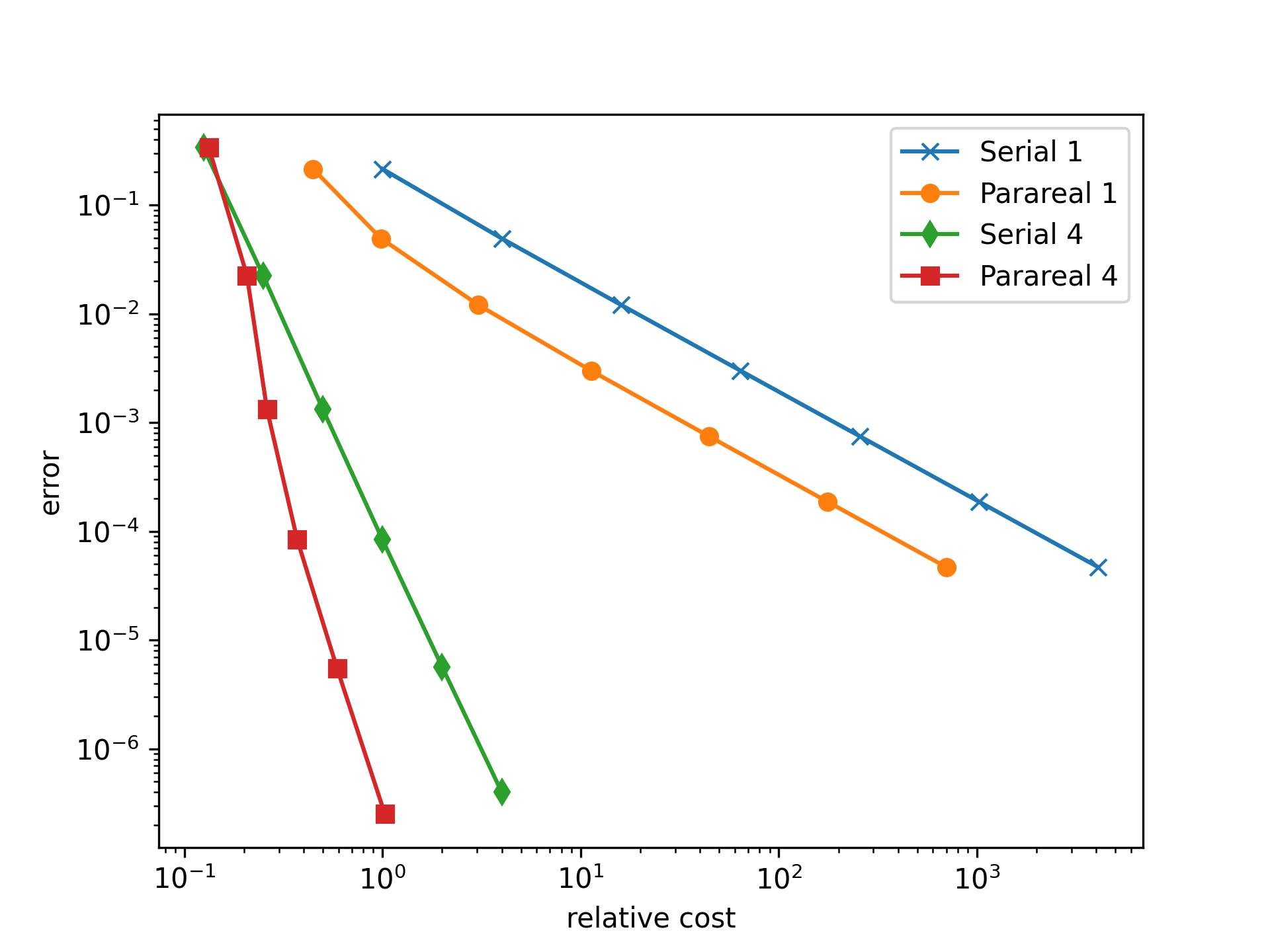}
     \end{subfigure}
    \caption{Comparison of serial and Parareal execution time for the K-S example using a first- and fourth-order ERK integrators. Note that the serial fourth-order integrator is always faster for a given accuracy than the parallel first-order method ({\it Way} 11)}
    \label{fig:KS}
\end{figure}

For this example we will again use Parareal applied to the Kuramoto-Siva\-shinsky Equation.  The K-S equation is a good choice 
to impress an audience because it gives rise to  chaotic temporal dynamics
(avoiding \textit{Way} 1). 
The equation in one dimension reads
\begin{align*}
    u_t = - u u_x - u_{xx} - u_{xxxx},
\end{align*}
which we solve on the spatial interval $x \in [0,32\pi]$ and temporal interval $t \in [0,40]$.  Since the fourth-derivative term is linear and stiff, we choose a first order exponential integrator in a spectral-space discretization where the 
linear operators diagonalize and hence the exponential of the operator is trivial to compute.  We use 512 points in space, and in this study will compare a serial first-order method with Parareal using the same first-order method in terms of cost per accuracy.  Using 32 time processors for all runs, we increase the number of steps for the fine Parareal propagator $\Fprop$ and hence the total number of time steps. The theoretical speedups (ignoring {\it Way} 9) are displayed on the left panel of Fig. \ref{fig:KS}.  One can see that the Parareal method provides speedup at all temporal resolutions up to a maximum of about 5.85 at the finest resolution (where $\alpha$ is the smallest).  So we have achieved meaningful speedup with a respectable efficiency for a problem with complex dynamics.  Best to stop the presentation here.

If we are a little more ambitious, we might replace our first-order integrator with the 4th-order exponential Runge-Kutta (ERK) method from \cite{Krogstad2005-kc}.  Now we need to be  more careful about {\it Way} 5 and hence won't be able to take nearly as many time steps.
In the right panel of Fig. \ref{fig:KS} we show the 1st-order and 4th-order results together.  The maximum theoretical speedup attained with the 4th-order method is only about 3.89 at the finest resolution,  which  is probably reason enough not to do the comparison.  But there is the additional irritation that at any accuracy level, the serial 4th-order method is significantly faster than the Parareal 1st-order method.  
\begin{way}
It is best to show speedup for first-order time integrators since they are a bit easier to inflate.  If you want to show speedup for higher-order methods as well, make it impossible to compare cost versus accuracy between first-order and higher-order methods.
\end{way}

\section{Parting thoughts}

The careful reader may have noticed that in all the examples above, a single PinT method is used for each {\it Way}.  
This brings  us finally to 
\begin{way}
 Never compare your own PinT method to a different PinT method.
\end{way}
The problem, as we have seen, is that assessing performance for a single PinT method is already not straightforward.
Comparing the performance of two  or more different methods makes matters even more difficult.
Although it has been often discussed within the PinT community, efforts to establish a set of benchmark test examples have, to date, made little head way.  
The performance of methods like PFASST and Parareal considered here are highly sensitive to the type of equation being solved, the type of spatial discretization being used, the accuracy desired, and the choice of problem and method parameters.  
In this study we purposely choose examples that lead to inflated reported speedups, and doing this required us to use our understanding of the methods and the equations chosen. 
Conversely, in most instances, a simple change in the experiment leads to much worse reported speedups.  
Different PinT approaches have strengths and weaknesses for different benchmark scenarios, hence establishing a set of benchmarks that the community would find  {\it fair} is a very non-trivial problem.  

Roughly, the ways we present can be grouped into three categories: ``choose your problem'' (ways 1--4), ``over-resolve'' (ways 5--8) and ``choose your performance measure'' (ways 9--11).  This classification is not perfect as some of the {\it Ways} overlap. 
Some of the dubious tricks presented here are intentionally obvious to detect, while others are more subtle. As in the original ``twelve ways'' article, and those it inspired, the examples are meant to be light-hearted. However, many of the \textit{Ways} have been (unintentionally) used when reporting numerical results, and the authors are not without guilt in this respect. 
Admitting that, we hope this article will be read the way we intended: as a demonstration of some of the many pitfalls one faces when assessing PinT performance and a reminder that considerable care is required to obtain truly meaningful results.

\section*{Acknowledgements}
The work of Minion  was supported by the U.S. Department of Energy, Office of Science,
Office of Advanced Scientific Computing Research, Applied Mathematics program
under contract number DE-AC02005CH11231. Part of the simulations were performed
using resources of the National Energy Research Scientific Computing Center
(NERSC), a DOE Office of Science User Facility supported by the Office of
Science of the U.S. Department of Energy under Contract No. DE-AC02-05CH11231.

\bibliographystyle{spmpsci}      
\bibliography{12ways,pint} 

\begin{thebibliography}{10}
\providecommand{\url}[1]{{#1}}
\providecommand{\urlprefix}{URL }
\expandafter\ifx\csname urlstyle\endcsname\relax
  \providecommand{\doi}[1]{DOI~\discretionary{}{}{}#1}\else
  \providecommand{\doi}{DOI~\discretionary{}{}{}\begingroup
  \urlstyle{rm}\Url}\fi

\bibitem{aktosun2007exact}
Aktosun, T., Demontis, F., Van Der~Mee, C.: Exact solutions to the focusing
  nonlinear {Schr{\"o}dinger} equation.
\newblock Inverse Problems \textbf{23}(5), 2171 (2007)

\bibitem{Ascher1997-er}
Ascher, U.M., Ruuth, S.J., Spiteri, R.J.: Implicit-explicit {Runge-Kutta}
  methods for time-dependent partial differential equations.
\newblock Appl. Numer. Math. \textbf{25}(2), 151--167 (1997)

\bibitem{Bailey1991-pr}
Bailey, D.H.: Twelve ways to fool the masses when giving performance results on
  parallel computers.
\newblock Supercomputing Review \textbf{4}(8) (1991)

\bibitem{Chawner2011-bj}
Chawner, J.: Revisiting ``twelve ways to fool the masses when describing mesh
  generation performance''.
\newblock
  \url{https://blog.pointwise.com/2011/05/23/revisiting-%e2%80%9ctwelve-ways-to-fool-the-masses-when-describing-mesh-generation-performance%e2%80%9d/}
  (2011).
\newblock Accessed: 2020-4-28

\bibitem{DongarraEtAl2014}
Dongarra, J., et~al.: {Applied Mathematics Research for Exascale Computing}.
\newblock Tech. Rep. LLNL-TR-651000, Lawrence Livermore National Laboratory
  (2014).
\newblock
  \urlprefix\url{{http://science.energy.gov/%7E/media/ascr/pdf/research/am/docs/EMWGreport.pdf}}

\bibitem{EmmettMinion2012}
Emmett, M., Minion, M.L.: {Toward an Efficient Parallel in Time Method for
  Partial Differential Equations}.
\newblock Communications in Applied Mathematics and Computational Science
  \textbf{7}, 105--132 (2012).
\newblock \doi{10.2140/camcos.2012.7.105}.
\newblock \urlprefix\url{http://dx.doi.org/10.2140/camcos.2012.7.105}

\bibitem{Gander2015_Review}
Gander, M.J.: {50 years of Time Parallel Time Integration}.
\newblock In: Multiple Shooting and Time Domain Decomposition. Springer (2015).
\newblock \doi{10.1007/978-3-319-23321-5_3}.
\newblock \urlprefix\url{http://dx.doi.org/10.1007/978-3-319-23321-5_3}

\bibitem{Gear1988}
Gear, C.W.: {Parallel methods for ordinary differential equations}.
\newblock CALCOLO \textbf{25}(1-2), 1--20 (1988).
\newblock \doi{10.1007/BF02575744}.
\newblock \urlprefix\url{http://dx.doi.org/10.1007/BF02575744}

\bibitem{299418}
{Globus}, A., {Raible}, E.: Fourteen ways to say nothing with scientific
  visualization.
\newblock Computer \textbf{27}(7), 86--88 (1994)

\bibitem{Gustafson1991-wi}
Gustafson, J.L.: Twelve ways to fool the masses when giving performance results
  on traditional vector computers.
\newblock \url{http://www.johngustafson.net/fun/fool.html} (1991).
\newblock Accessed: 2020-4-28

\bibitem{Hoefler2018-wb}
Hoefler, T.: Twelve ways to fool the masses when reporting performance of deep
  learning workloads! (not to be taken too seriously) (2018).
\newblock ArXiv, 1802.09941

\bibitem{Krogstad2005-kc}
Krogstad, S.: Generalized integrating factor methods for stiff {PDEs}.
\newblock J. Comput. Phys. \textbf{203}(1), 72--88 (2005)

\bibitem{LionsEtAl2001}
Lions, J.L., Maday, Y., Turinici, G.: {A "parareal" in time discretization of
  {PDE}'s}.
\newblock Comptes Rendus de l'Académie des Sciences - Series I - Mathematics
  \textbf{332}, 661--668 (2001).
\newblock \doi{10.1016/S0764-4442(00)01793-6}.
\newblock \urlprefix\url{http://dx.doi.org/10.1016/S0764-4442(00)01793-6}

\bibitem{Minhas2019-rm}
Minhas, F., Asif, A., Ben-Hur, A.: Ten ways to fool the masses with machine
  learning (2019).
\newblock ArXiv, 1901.01686

\bibitem{Nievergelt1964}
Nievergelt, J.: {Parallel methods for integrating ordinary differential
  equations}.
\newblock Commun. ACM \textbf{7}(12), 731--733 (1964).
\newblock \doi{10.1145/355588.365137}.
\newblock \urlprefix\url{http://dx.doi.org/10.1145/355588.365137}

\bibitem{ong2019applications}
Ong, B.W., Schroder, J.B.: Applications of time parallelization.
\newblock Computing and Visualization in Science \textbf{648} (2019)

\bibitem{Pakin2011-ec}
Pakin, S.: Ten ways to fool the masses when giving performance results on
  {GPUs}.
\newblock HPCwire, December \textbf{13} (2011)

\bibitem{Tautges2004-ry}
Tautges, T.J., White, D.R., Leland, R.W.: Twelve ways to fool the masses when
  describing mesh generation performance.
\newblock IMR/PINRO Joint Rep. Ser. pp. 181--190 (2004)

\end{thebibliography}

\section*{Appendix 1} 
\appendix
The value of $\sigma$ in the first example is 0.02.

\end{document}